\documentclass[aos,preprint]{imsart}
\RequirePackage{amsthm,amsmath,amsfonts,amssymb}
\RequirePackage[numbers]{natbib}
\RequirePackage{mathtools}
\RequirePackage{graphicx}
\RequirePackage[colorlinks,citecolor=blue,urlcolor=blue]{hyperref}

\startlocaldefs
\numberwithin{equation}{section}
\newtheorem{theorem}{Theorem}[section]
\newtheorem{proposition}[theorem]{Proposition}
\newtheorem{lemma}[theorem]{Lemma}
\newtheorem{corollary}[theorem]{Corollary}
\theoremstyle{definition}

\newcommand{\R}{\mathbb R}
\newcommand{\E}{\mathbb E}

\newcommand{\dd}{\,\mathrm d}
\newcommand{\one}{\mathbf 1}
\newcommand{\norm}[1]{\left\lVert #1\right\rVert}

\usepackage{xcolor}

\DeclareMathOperator{\Var}{Var}
\DeclareMathOperator{\Cov}{Cov}
\DeclareMathOperator{\diag}{diag}
\endlocaldefs

\begin{document}
\begin{frontmatter}
\title{Markov potential priors for prediction in Gamma models}
\runtitle{Shrinkage prediction in Gamma models}
\runauthor{X. Li and F. Komaki}
\pdfsubject{Author manuscript prepared for The Annals of Statistics}
\begin{aug}
\author{\fnms{Xiao}~\snm{Li}\ead[label=e1]{lixiaoms@163.com}}
\author[B,C]{\fnms{Fumiyasu}~\snm{Komaki}\ead[label=e2]{komaki@mist.i.u-tokyo.ac.jp}}
\address{\printead{e1} \printead{e2}}
\address[B]{Department of Mathematical Informatics,
The University of Tokyo, Tokyo, Japan}
\address[C]{Riken Center for Brain Science, Saitama, Japan}
\end{aug}

\begin{abstract}
We consider Bayesian prediction in independent Gamma scale models with
arbitrary known positive shape parameters. When the dimension is at least
three, we construct a shrinkage prior toward any prespecified positive scale
vector. The resulting predictive density strictly dominates the Jeffreys
predictive density under Kullback--Leibler loss and is minimax.
The prior is obtained by weighting the Jeffreys prior with a Markov
potential. The transition is determined by the beta--gamma decomposition.
The domination result is nonasymptotic.
\end{abstract}

\begin{keyword}[class=MSC]
\kwdgroup[type=primary]{\kwd{62C15}\kwd{62F15}}
\kwdgroup[type=secondary]{\kwd{62C10}\kwd{62C20}}
\end{keyword}
\begin{keyword}
\kwd{Bayesian prediction}
\kwd{Gamma scale model}
\kwd{Kullback--Leibler loss}
\kwd{Markov potential}
\kwd{Shrinkage prior}
\kwd{Minimax prediction}
\end{keyword}
\end{frontmatter}

\section{Introduction}\label{sec:introduction}

We consider prediction in independent Gamma models with known positive
shape parameters and an unknown scale in each coordinate. The observed
and future shape vectors are fixed and may differ in every coordinate.
Predictive densities are compared by the expected Kullback--Leibler
divergence from the true future distribution. We construct shrinkage
priors whose predictive densities strictly dominate the product Jeffreys
rule and are minimax when the dimension is at least three.

In a multiparameter prediction problem,
a predictive density obtained from an invariant prior can be improved by
introducing dependence across the parameters through the prior.
For normal observations, \citet{komaki2001} established a predictive
shrinkage phenomenon, and \citet{georgeliangxu2006} developed improved
minimax predictive densities under Kullback--Leibler loss.
For independent Poisson observations, \citet{komaki2004} obtained
finite-sample comparisons of Bayesian predictive densities.
\citet{likomaki2026} give superharmonic sufficient conditions for
improving Jeffreys predictive distributions in independent Poisson
models, including unequal observation and prediction durations.
These results motivate the search for priors adapted to the sampling model
and to the distinction between the observed and future experiments.

Classical work also treats inadmissibility in simultaneous estimation of
Gamma scale parameters. \citet{berger1980} developed
differential-inequality methods for constructing improvements, and
\citet{dasgupta1986} studied simultaneous estimation under weighted
quadratic losses. Those results concern vector-valued parameter
estimators. The decision object here is a predictive density and the loss
is Kullback--Leibler divergence, so the relevant risk comparison is
established directly for that decision problem.

For Gamma models, \citet{lmoudden2017} studied predictive density
estimation with scale constraints and identified the unrestricted scalar
minimax benchmark. \citet{komaki2026} established admissibility of the
scalar Jeffreys predictive density for arbitrary positive observed and
future shape parameters.

Potential theory provides a useful source of such priors.
\citet{komaki2006} related shrinkage priors and Green functions
on statistical manifolds in an asymptotic predictive framework.
A useful construction is the potential of a Markov transition
operator, obtained by summing its iterates applied to a
nonnegative function. Such potentials and their Poisson
equations are classical; see
\citet[Section~4.3]{lawlerlimic2010}.
In the present Gamma model, the beta--gamma decomposition
naturally determines a Markov transition operator $Q_{ab}$,
where $a$ and $b$ denote the observed and future shape vectors.
Its increments are differences of independent negative
logarithms of Beta variables.
For $d\ge3$, the potential of a strictly positive Gaussian
function $f$ satisfies
\[
 w_{ab}=\sum_{n=0}^{\infty}Q_{ab}^n f,
 \qquad
 w_{ab}-Q_{ab}w_{ab}=f>0.
\]
The second identity is the Poisson equation for $Q_{ab}$,
and its strict positivity provides the key condition
for finite-sample predictive risk improvement.
The operator and the potential are defined precisely below.

For any prespecified $\beta_0\in(0,\infty)^d$, our prior reweights
the flat Jeffreys prior in logarithmic scale coordinates.
Its weight is symmetric about $\theta_0=\log\beta_0$ and
has its strict global maximum there.
It depends on the known shape parameters and the target $\beta_0$.
For fixed shape parameters, its far tail has the classical harmonic shrinkage
order $\norm{\theta-\theta_0}^{2-d}$, with an explicit anisotropic
leading term.
Dependence of a prior on a fixed observation and prediction design also
occurs in the Poisson-process setting of \citet{komaki2015}.

The proof proceeds in two main steps. First, the beta--gamma
coupling and conditional Jensen's inequality give a sufficient
condition for predictive risk improvement. Second, Fourier analysis
constructs a bounded positive weight satisfying this condition,
while logarithmic moment bounds justify the required expectations.

The main result and minimax benchmark are stated in
Section~\ref{sec:model-main}.
Section~\ref{sec:criterion} gives the exact coupling, predictive formula,
and domination criterion. Section~\ref{sec:potential} constructs the
potential and proves the main result, including a quantitative refinement.
Section~\ref{sec:example} treats a three-dimensional exponential example
and illustrates its risk along a common-scale slice. The appendices give
the distributional calculations, kernel estimates, pointwise Fourier
calculation, integrability arguments, minimax lower bound, and numerical
method.

\section{Model and main result}\label{sec:model-main}

Fix an integer $d\ge1$ and known vectors $a,b\in(0,\infty)^d$,
and write $c=a+b$. The parameter is
$\beta\in\mathcal B=(0,\infty)^d$.
Given $\beta$, all the variables $X_i,Y_i$ are independent, with
\begin{equation*}\label{eq:model}
 X_i\sim\operatorname{Gamma}(a_i,\text{scale }\beta_i),
 \qquad
 Y_i\sim\operatorname{Gamma}(b_i,\text{scale }\beta_i).
\end{equation*}
Thus, for any positive shape vector $r$, the density with respect to
Lebesgue measure on $(0,\infty)^d$ is
\begin{equation*}\label{eq:gamma-density}
 p_r(z\mid\beta)=
 \prod_{i=1}^d\frac{z_i^{r_i-1}e^{-z_i/\beta_i}}
                   {\Gamma(r_i)\beta_i^{r_i}},
 \qquad z\in(0,\infty)^d.
\end{equation*}
All vector logarithms, exponentials, products, and quotients below are
coordinatewise unless indicated otherwise. We use the Euclidean norm
$\norm{\cdot}$, set $\theta=\log\beta$, and keep the coordinates of
$\beta$ free throughout.

The Kullback--Leibler risk of a predictive density $q(y\mid x)$ is
\begin{equation}\label{eq:risk}
 R(\beta,q)=\E_{X,Y\mid\beta}
       \log\frac{p_b(Y\mid\beta)}{q(Y\mid X)}.
\end{equation}
Write $\mathcal Q$ for the class of all jointly measurable predictive
densities normalized in $y$ for every positive $x$.

The Jeffreys prior is
$\pi_J(\beta)\dd\beta=\prod_i\beta_i^{-1}\dd\beta=\dd\theta$.
Its Bayesian predictive density is
\begin{equation}\label{eq:qJ}
 q_J(y\mid x)=\prod_{i=1}^d
 \frac{\Gamma(c_i)}{\Gamma(a_i)\Gamma(b_i)}
 \frac{x_i^{a_i}y_i^{b_i-1}}{(x_i+y_i)^{c_i}}.
\end{equation}
Write $\psi=(\log\Gamma)'$. The constant risk of \eqref{eq:qJ} is
\begin{equation}\label{eq:CJ}
 C_J=\sum_{i=1}^d
 \bigl[\log\Gamma(a_i)-\log\Gamma(c_i)-b_i
       -a_i\psi(a_i)+c_i\psi(c_i)\bigr].
\end{equation}
The risk calculation is given in Lemma~\ref{lem:original-losses}.

For a positive weight function $w$ on $\R^d$, consider the prior
\begin{equation}\label{eq:prior}
 \pi_w(\beta)\dd\beta
 =w(\log\beta)\prod_{i=1}^d\beta_i^{-1}\dd\beta
 =w(\theta)\dd\theta.
\end{equation}
Let $G_r$ have independent $\operatorname{Gamma}(r_i,1)$ coordinates,
write
\begin{equation}\label{eq:LG}
 L_r=\log G_r,
\end{equation}
and define
\begin{equation}\label{eq:n-definition}
 n_r(u)=\E w(u-L_r).
\end{equation}
When $w$ is positive and bounded, posterior propriety and predictive
normalization follow from Lemma~\ref{lem:posterior}. The corresponding
Bayesian predictive density is
\begin{equation}\label{eq:qw}
 q_w(y\mid x)=q_J(y\mid x)
       \frac{n_c(\log(x+y))}{n_a(\log x)}.
\end{equation}

The prior in our main result is specified by a probability transition.

Take independent
$B_{a_i b_i}\sim\operatorname{Beta}(a_i,b_i)$,
and define the random vector
\begin{equation}\label{eq:LB}
 D_{ab}=(D_{a_i b_i})_{i=1}^d,
 \qquad D_{a_i b_i}=-\log B_{a_i b_i}.
\end{equation}
Let $D'_{ab}$ be an independent copy of $D_{ab}$.
For any bounded measurable function $h:\R^d\to\R$, define
\begin{equation}\label{eq:Q}
 Q_{ab}h(u)=\E h(u+D'_{ab}-D_{ab}).
\end{equation}
Fix a prespecified scale target $\beta_0\in(0,\infty)^d$ and set
$\theta_0=\log\beta_0$. Let
\begin{equation}\label{eq:potential}
 f(u)=f_{\theta_0}(u)=e^{-\norm{u-\theta_0}^2},
 \qquad
 w_{ab}(u)=\sum_{n=0}^{\infty}Q_{ab}^nf(u),
\end{equation}
where $Q_{ab}^0$ is the identity operator.
We call $w_{ab}$ the Markov potential, which defines
a shrinkage prior density through \eqref{eq:prior}.

\begin{theorem}\label{thm:main}
Suppose $d\ge3$, $a_i,b_i>0$, and
$\beta_0\in(0,\infty)^d$ is fixed.
Set $\theta_0=\log\beta_0$.
Then the following statements hold.
\begin{enumerate}
\item
The series in \eqref{eq:potential} converges uniformly on $\R^d$.
Its sum is continuous, and
\begin{equation*}\label{eq:main-bounds}
 M:=\sum_{n=0}^{\infty}\norm{Q_{ab}^nf}_\infty<\infty,
 \qquad 0<f\le w_{ab}\le M,
 \qquad Q_{ab}w_{ab}=w_{ab}-f.
\end{equation*}
The prior measure \eqref{eq:prior} with $w=w_{ab}$
is improper.

\item
For every $x\in(0,\infty)^d$, the posterior under
\eqref{eq:prior} with $w=w_{ab}$ is proper.
The corresponding predictive density $q_{w_{ab}}$ is given by
\eqref{eq:qw}, is positive and finite for every positive
$(x,y)$, and integrates to one in $y$ for every positive $x$.

\item
For every $\beta\in(0,\infty)^d$ and
$q\in\{q_J,q_{w_{ab}}\}$,
\begin{equation}\label{eq:main-absolute-loss}
 \E_{X,Y\mid\beta}
 \left|\log\frac{p_b(Y\mid\beta)}{q(Y\mid X)}\right|
 <\infty,
\end{equation}
and
\begin{equation}\label{eq:main-domination}
 0\le R(\beta,q_{w_{ab}})<R(\beta,q_J)<\infty.
\end{equation}
\item
Both $q_J$ and $q_{w_{ab}}$ are minimax over $\mathcal Q$, with
\begin{equation}\label{eq:main-minimax}
 \inf_{q\in\mathcal Q}\sup_{\beta\in\mathcal B}R(\beta,q)
 =\sup_{\beta\in\mathcal B}R(\beta,q_{w_{ab}})
 =C_J,
\end{equation}
where $C_J$ is defined by \eqref{eq:CJ}.
\end{enumerate}
Consequently, $q_J$ is inadmissible in the class of all
predictive densities.
\end{theorem}

The shape parameters may differ across coordinates and between the
observation and future samples. In particular, the shape parameters
need not be integers and may be smaller than one.
The quantitative bound in Corollary~\ref{cor:quantitative}
is a refinement of strict domination; the basic comparison is driven
by the identity $w_{ab}-Q_{ab}w_{ab}=f$.
Relative to the flat measure $\dd\theta$, the prior gives its
largest weight to $\theta_0=\log\beta_0$; the precise symmetry
and strict-maximum properties are given in
\eqref{eq:centered-weight}.

\subsection{The minimax benchmark}

Jeffreys minimaxity is a consequence of the general location-model result
of \citet[Appendix~B, Theorem~1$'$]{liangbarron2004},
applied after the coordinatewise logarithmic transformation.
The observed and future location noises are $\log G_a$ and
$\log G_b$ and may have different distributions.
The scalar Gamma calculation is also given by
\citet[Corollary~3.1]{lmoudden2017}.
For completeness, we give a direct proof of
the lower bound over the full class $\mathcal Q$.

\begin{proposition}\label{prop:minimax}
For every positive integer $d$ and positive shape vectors $a,b$,
\begin{equation*}\label{eq:minimax-value}
 \inf_{q\in\mathcal Q}\sup_{\beta\in\mathcal B}R(\beta,q)=C_J,
\end{equation*}
where $C_J$ is defined by \eqref{eq:CJ}.
Thus, the Jeffreys predictive density $q_J$ is minimax.
Every predictive rule that dominates $q_J$ is also minimax.
\end{proposition}

\begin{proof}
Fix $\tau_i>0$ and give the rates $\lambda_i=1/\beta_i$
independent $\operatorname{Gamma}(\varepsilon,\text{rate }\tau_i)$
priors, where $\varepsilon>0$.
Their joint prior is proper. Its Bayes risk is denoted by
$r_\varepsilon$.
Appendix~\ref{app:minimax} calculates $r_\varepsilon$
from the logarithmic loss and proves
$r_\varepsilon\to C_J$ as $\varepsilon\downarrow0$.
The Bayes-risk decomposition in Appendix~\ref{app:minimax} gives
$\sup_{\beta\in\mathcal B}R(\beta,q)\ge r_\varepsilon$
for every $q\in\mathcal Q$ and every $\varepsilon>0$.
Letting $\varepsilon\downarrow0$ yields the minimax lower bound.

The constant risk $C_J$ of $q_J$ supplies the matching upper bound.
If $q$ dominates $q_J$, then
$\sup_\beta R(\beta,q)\le C_J$;
the same lower bound forces equality.
\end{proof}

\section{An exact coupling and a domination criterion}
\label{sec:criterion}

\subsection{The beta--gamma coupling}

The beta--gamma decomposition gives independent Gamma and Beta variables
$T_i=U_i+V_i$ and $B_i=U_i/T_i$, with shape parameters $c_i$ and
$(a_i,b_i)$, respectively. Here $U_i,V_i$ are independent unit-scale
Gamma variables with shape parameters $a_i,b_i$.

Applying the log transformation in each coordinate gives an exact coupling
\begin{equation}\label{eq:coupling}
 L_a=L_c-D_{ab},\qquad L_c\ \hbox{independent of }D_{ab},
\end{equation}
where $L_r$ and $D_{ab}$ are defined by
\eqref{eq:LG} and \eqref{eq:LB}, respectively.
It holds for every positive pair of shape
parameter vectors $a,b$.

Define the probability convolution operators
\begin{equation*}\label{eq:operators}
 A_rh(u)=\E h(u-L_r),\quad
 P_{ab+}h(u)=\E h(u+D_{ab}),\quad
 P_{ab-}h(u)=\E h(u-D_{ab}).
\end{equation*}
They are defined for bounded measurable functions and for nonnegative
measurable functions with extended values.
These convolution operators commute. Together with the coupling
\eqref{eq:coupling}, this gives
\begin{equation}\label{eq:operator-identities}
 A_a=P_{ab+}A_c,\qquad
 Q_{ab}=P_{ab-}P_{ab+}=P_{ab+}P_{ab-},
 \qquad A_rQ_{ab}=Q_{ab}A_r.
\end{equation}
The sign in $A_a=P_{ab+}A_c$ follows from
$u-L_a=u-L_c+D_{ab}$.

\subsection{Posterior and predictive distributions}

\begin{lemma}\label{lem:posterior}
Let $w:\R^d\to(0,\infty)$ be bounded and measurable.
Then \eqref{eq:prior} has a proper posterior for every positive $x$,
and its Bayesian predictive density is \eqref{eq:qw}.
This density is positive and finite at every
positive $(x,y)$, and integrates to one in $y$ for every positive $x$.
\end{lemma}

The proof in Appendix~\ref{app:distributional} gives the explicit
marginal and posterior normalizing constants.

\subsection{A sufficient condition for domination}

For the following proposition, the dimension may be any positive integer.

\begin{proposition}\label{prop:criterion}
Let $w:\R^d\to(0,\infty)$ be bounded and measurable,
and recall that $n_r=A_rw$.
Suppose
\begin{align}
 &\E|\log n_r(\theta+L_r)|<\infty
   \quad\text{for every }\theta\in\R^d,\quad r\in\{a,c\},
   \label{eq:criterion-outer}\\
 &P_{ab-}|\log n_a|(u)<\infty
   \quad\text{for every }u\in\R^d.
   \label{eq:criterion-inner}
\end{align}
If $Q_{ab}w\le w$ pointwise, then $q_w$ is a predictive density,
and $R(\beta,q_w)\le R(\beta,q_J)$ for every $\beta$.
If $Q_{ab}w<w$ at every point, the risk inequality is strict
at every $\beta$.
\end{proposition}

\begin{proof}
Lemma~\ref{lem:posterior}
establishes the existence of the predictive density.
Lemma~\ref{lem:original-losses}
and the integrability arguments in
Appendix~\ref{app:endpoint-regularity} show that the
two logarithmic losses are absolutely integrable
and have nonnegative risks.

From \eqref{eq:qw}, we may subtract the risks and use the sampling
marginals $\log X\overset{\mathcal D}=\theta+L_a$ and
$\log(X+Y)\overset{\mathcal D}=\theta+L_c$:
\begin{equation}\label{eq:risk-marginals}
 \Delta(\beta)\coloneqq R(\beta,q_J)-R(\beta,q_w)
 =\E\log n_c(\theta+L_c)-\E\log n_a(\theta+L_a).
\end{equation}

The first identity in \eqref{eq:operator-identities}
gives $n_a=P_{ab+}n_c$.
The independent coupling \eqref{eq:coupling} and
the justified Fubini step in Appendix~\ref{app:endpoint-regularity}
yield the identity
\begin{equation}\label{eq:endpoint-identity}
 \Delta(\beta)
 =\E_{L_c}\{\log n_c-P_{ab-}\log(P_{ab+}n_c)\}(\theta+L_c),
\end{equation}
where the relation
$\theta+L_a=\theta+L_c-D_{ab}$
is used.

Fix $u$. The random variable
$Z_u=(P_{ab+}n_c)(u-D_{ab})=n_a(u-D_{ab})$
is positive and bounded, and \eqref{eq:criterion-inner} gives
$\E|\log Z_u|<\infty$. Thus Jensen's inequality applies to its logarithm:
\begin{equation*}\label{eq:endpoint-Jensen}
 P_{ab-}\log(P_{ab+}n_c)(u)
 \le\log\{P_{ab-}P_{ab+}n_c(u)\}
 =\log Q_{ab}n_c(u).
\end{equation*}

The assumed condition $Q_{ab}w\le w$ gives
$Q_{ab}n_c=A_cQ_{ab}w\le n_c$.
Both $n_c$ and $Q_{ab}n_c$ are everywhere positive and finite. Hence
\begin{equation}\label{eq:pointwise-risk-improvement}
 A(u)\coloneqq \log n_c(u)
 -P_{ab-}\log(P_{ab+}n_c)(u)
 \ge\log\frac{n_c(u)}{Q_{ab}n_c(u)}\ge0.
\end{equation}

By \eqref{eq:A-integrable},
$\E|A(\theta+L_c)|<\infty$.
Taking expectations in \eqref{eq:pointwise-risk-improvement} and using
\eqref{eq:endpoint-identity} proves weak domination.

If $w-Q_{ab}w>0$ everywhere, then
$n_c-Q_{ab}n_c=A_c(w-Q_{ab}w)>0$ everywhere,
so $A(u)>0$ at every $u$.
Since $A(\theta+L_c)>0$ almost surely and is integrable, we have
\[
 \Delta(\beta)=\E_{L_c}[A(\theta+L_c)]>0.
\]
Thus, strict domination holds for every $\theta$.
\end{proof}

Strict domination follows from the strict superharmonicity
of $w$, without requiring a strict version of
Jensen's inequality.

\section{Construction of the Markov potential}\label{sec:potential}

We now verify the sufficient conditions for the weight
\eqref{eq:potential}. Let $W_{ab}=D_{ab}'-D_{ab}$, and write
\begin{equation}\label{eq:characteristic}
 \phi_{ab}(t)=\E e^{it\cdot W_{ab}}
  =\left|\E e^{it\cdot D_{ab}}\right|^2\in[0,1].
\end{equation}
The covariance of $W_{ab}$ is $2K$, where
\begin{equation}\label{eq:K}
 K=\diag(K_1,\ldots,K_d),\qquad
 K_i=\psi_1(a_i)-\psi_1(c_i)>0,
\end{equation}
and $\psi=(\log\Gamma)'$, $\psi_1=\psi'$.
Appendix~\ref{app:kernel} proves these statements directly for all
positive shapes, together with the following bounds: there exist
$\delta,\kappa>0$ and $0<\rho<1$, depending only on $a,b$, such that
\begin{equation}\label{eq:frequency-bounds}
 \phi_{ab}(t)\le e^{-\kappa\norm{t}^2}
 \quad(\norm{t}\le\delta),
 \qquad
 \sup_{\norm{t}\ge\delta}\phi_{ab}(t)\le\rho.
\end{equation}

\begin{lemma}\label{lem:iterate-fourier}
For every integer $n\ge0$ and every $u\in\R^d$,
\begin{equation}\label{eq:Qn-pointwise-fourier}
 Q_{ab}^nf(u)=(2\pi)^{-d}\int_{\R^d}
    e^{it\cdot(u-\theta_0)}\pi^{d/2}e^{-\norm{t}^2/4}
       \phi_{ab}(t)^n\dd t.
\end{equation}
The integral is absolutely convergent.
\end{lemma}

The proof in Appendix~\ref{app:Fourier} uses Gaussian inversion and
Fubini's theorem and applies in every positive dimension.

\begin{proposition}\label{prop:potential}
For $n\ge1$, the Gaussian function in \eqref{eq:potential} satisfies
\begin{equation}\label{eq:iterate-bound}
 \norm{Q_{ab}^nf}_\infty
 \le2^{-d}\kappa^{-d/2}n^{-d/2}+\rho^n.
\end{equation}
If $d\ge3$, its potential is continuous and obeys
\begin{equation}\label{eq:M-bound}
 1\le M:=\sum_{n=0}^{\infty}\norm{Q_{ab}^nf}_\infty
 \le1+2^{-d}\kappa^{-d/2}\frac{d}{d-2}
       +\frac{\rho}{1-\rho}<\infty.
\end{equation}
Furthermore,
\begin{equation}\label{eq:potential-defect}
 0<f\le w_{ab}\le M,
 \qquad Q_{ab}w_{ab}=w_{ab}-f,
 \qquad \int_{\R^d}w_{ab}(u)\dd u=\infty.
\end{equation}
\end{proposition}

\begin{proof}
Lemma~\ref{lem:iterate-fourier} and \eqref{eq:frequency-bounds} give,
for $n\ge1$,
\begin{align*}
 \norm{Q_{ab}^nf}_\infty
 &\le(2\pi)^{-d}\int_{\R^d}
          \pi^{d/2}e^{-\norm{t}^2/4}\phi_{ab}(t)^n\dd t\\
 &\le(2\pi)^{-d}\pi^{d/2}
          \int_{\norm{t}\le\delta}e^{-n\kappa\norm{t}^2}\dd t
       +\rho^n(2\pi)^{-d}
          \int_{\norm{t}>\delta}\pi^{d/2}e^{-\norm{t}^2/4}\dd t\\
 &\le2^{-d}\kappa^{-d/2}n^{-d/2}+\rho^n.
\end{align*}
The last step enlarges both regions to $\R^d$ and evaluates the
Gaussian integrals. This proves \eqref{eq:iterate-bound}.

For the remaining claims, assume $d\ge3$.
Since $f$ is bounded and continuous, each iterate
$Q_{ab}^nf$ is continuous by the bounded convergence theorem.
Since
\[
 \sum_{n=1}^{\infty}n^{-d/2}
 \le1+\int_1^\infty s^{-d/2}\dd s=\frac{d}{d-2},
\]
the preceding bounds on the iterates, together with
$\norm{f}_\infty=1$, give \eqref{eq:M-bound}.
The Weierstrass M-test then shows that the series
defining $w_{ab}$ converges uniformly, and its sum is continuous.
Since all terms are nonnegative, we have
$f\le w_{ab}\le M$.

Each of the two convolutions defining $Q_{ab}$ is an average
with respect to a probability measure and therefore cannot
increase the supremum norm.
Thus, $Q_{ab}$ can be applied term by term to the uniformly
convergent series, giving
\[
 Q_{ab}w_{ab}
 =\sum_{n\ge1}Q_{ab}^nf=w_{ab}-f.
\]

Finally, since convolution with a probability measure
preserves integrals, we have
\begin{equation*}\label{eq:mass-preserving}
 \int_{\R^d}Q_{ab}^nf(u)\dd u
 =\int_{\R^d}f(u)\dd u=\pi^{d/2}
 \qquad(n\ge0).
\end{equation*}
Consequently,
\[
 \int_{\R^d}w_{ab}(u)\dd u
 =\sum_{n=0}^{\infty}\int_{\R^d}Q_{ab}^nf(u)\dd u
 =\infty.
\]
\end{proof}

For $d\ge3$, Lemma~\ref{lem:iterate-fourier} and summation
of the geometric series in $\phi_{ab}(t)$ give the Fourier
integral representation
\begin{equation}\label{eq:potential-fourier}
 w_{ab}(u)=(2\pi)^{-d}\int_{\R^d}
       \cos\{t\cdot(u-\theta_0)\}
       \frac{\pi^{d/2}e^{-\norm{t}^2/4}}
            {1-\phi_{ab}(t)}\dd t.
\end{equation}
The bounds established in Appendix~\ref{app:Fourier}
ensure absolute convergence and justify the interchange
of summation and integration.

The Fourier representation \eqref{eq:potential-fourier} also shows that
\begin{equation}\label{eq:centered-weight}
 w_{ab}(\theta_0+v)=w_{ab}(\theta_0-v)<w_{ab}(\theta_0)
 \qquad(v\ne0).
\end{equation}
Indeed, symmetry follows from the cosine representation. Moreover,
for $v\ne0$,
\[
 w_{ab}(\theta_0)-w_{ab}(\theta_0+v)
 =(2\pi)^{-d}\int_{\R^d}
 [1-\cos(t\cdot v)]
 \frac{\pi^{d/2}e^{-\norm{t}^2/4}}
      {1-\phi_{ab}(t)}\dd t>0,
\]
since the integrand is positive almost everywhere.
Thus, in logarithmic scale coordinates, the weight relative to
the flat Jeffreys prior attains its unique maximum at
$\theta_0=\log\beta_0$.

\begin{proposition}\label{prop:potential-tail}
For $d\ge3$ and every fixed pair of positive shape vectors $a,b$, the
potential in \eqref{eq:potential} satisfies
\begin{equation}\label{eq:potential-tail}
 w_{ab}(\theta_0+v)\sim
 \frac{\Gamma(d/2-1)}{4\sqrt{\det K}}\,
       (v^{\mathsf T}K^{-1}v)^{(2-d)/2}
 \qquad\text{as }\norm{v}\longrightarrow\infty.
\end{equation}
The equivalence is uniform over unit directions for each fixed $K$.
In particular, $w_{ab}(\theta_0+v)\asymp\norm{v}^{2-d}$:
for all sufficiently large $\norm{v}$, their ratio is bounded above
and below by positive constants depending on the fixed shapes.
\end{proposition}

The proof in Appendix~\ref{app:tail} separates the quadratic
singularity of the Fourier multiplier from a faster-decaying
remainder. Formula~\eqref{eq:potential-tail} describes the weight
relative to $\dd\theta$. The prior density relative to $\dd\beta$
also includes the Jacobian factor $\prod_i\beta_i^{-1}$ in
\eqref{eq:prior}.

For $n_r=A_rw_{ab}$ and $h_r=A_rf$, the commutativity
of $A_r$ and $Q_{ab}$ gives
\begin{equation}\label{eq:smoothed-defect}
 Q_{ab}n_r=n_r-h_r,
 \qquad 0<h_r(u)<n_r(u)\le M.
\end{equation}
The inequalities follow from $h_r>0$, $Q_{ab}n_r>0$,
and $n_r\le M$.
This identity will yield a strict risk comparison once the
necessary logarithmic integrability is established.

\begin{lemma}\label{lem:log-bounds}
Let $w_{ab}$ be the potential in \eqref{eq:potential} with $d\ge3$,
and let $n_r=A_rw_{ab}$.
For a positive shape vector $r$, put
\begin{equation}\label{eq:log-moments}
 \mu_r=(\psi(r_i))_{i=1}^d=\E L_r,
 \qquad V_r=\sum_{i=1}^d\psi_1(r_i)
              =\E\norm{L_r-\mu_r}^2.
\end{equation}
Then $V_r<\infty$, and, for every $u,\theta\in\R^d$,
\begin{align}
 n_r(u)&\ge
      e^{-\norm{u-\theta_0-\mu_r}^2-V_r},
      \label{eq:smoothing-lower}\\
 |\log n_r(u)|
 &\le\log M+\norm{u-\theta_0-\mu_r}^2+V_r,
      \label{eq:log-envelope}\\
 \E|\log n_r(\theta+L_r)|
 &\le\log M+\norm{\theta-\theta_0}^2+2V_r.
      \label{eq:outer-log-bound}
\end{align}
In addition,
\begin{equation}\label{eq:inner-log-bound}
 P_{ab-}|\log n_a|(u)
 \le\log M+\E\norm{u-D_{ab}-\theta_0-\mu_a}^2
       +V_a<\infty.
\end{equation}
\end{lemma}

The key lower bound \eqref{eq:smoothing-lower} follows
from $w_{ab}\ge f$ and Jensen's inequality.
The remaining bounds follow from $n_r\le M$ and the
moment identities in \eqref{eq:log-moments}.
Details are given in Appendix~\ref{app:integrability}.
These bounds ensure the required logarithmic integrability
even though the potential need not be bounded away from zero.

\begin{proof}[Proof of Theorem~\ref{thm:main}]
Proposition~\ref{prop:potential} establishes uniform convergence,
continuity, the bound $M<\infty$, and the identity
$Q_{ab}w_{ab}=w_{ab}-f$.
It also shows that the prior \eqref{eq:prior} with $w=w_{ab}$
has infinite total mass and is therefore improper.

Lemma~\ref{lem:posterior} gives posterior propriety and a valid
predictive density.
Equations~\eqref{eq:outer-log-bound} and
\eqref{eq:inner-log-bound} verify both moment conditions of
Proposition~\ref{prop:criterion} for every $u$ and $\theta$.
Since $w_{ab}-Q_{ab}w_{ab}=f>0$,
Proposition~\ref{prop:criterion} gives
\eqref{eq:main-domination}.
The absolute integrability asserted in
\eqref{eq:main-absolute-loss} follows from
Lemma~\ref{lem:original-losses}.
The competing rule belongs to the full action class, so this proves
inadmissibility of $q_J$ in that class.

Proposition~\ref{prop:minimax} identifies the minimax value as $C_J$.
The pointwise bound $R(\beta,q_{w_{ab}})\le C_J$ and the minimax lower
bound then imply
$\sup_\beta R(\beta,q_{w_{ab}})=C_J$, proving
\eqref{eq:main-minimax}.

Finally, translations commute with $Q_{ab}$.
Each iterate of $f_{\theta_0}$ is a translate of the
corresponding iterate centered at zero.
Hence its supremum norm is independent of $\theta_0$,
and so is $M$.
\end{proof}

\subsection{A quantitative refinement}

The defect identity also gives a simple lower bound for the risk
improvement. The constant $M$ in the result below can be replaced by
any upper bound for it, such as the right side of \eqref{eq:M-bound}.

\begin{corollary}\label{cor:quantitative}
Under the conditions of Theorem~\ref{thm:main},
\begin{equation}\label{eq:quantitative-bound}
 \Delta(\beta)\ge M^{-1}
   \exp\left\{-\norm{\log\beta-\theta_0}^2
                    -2\sum_{i=1}^d\psi_1(c_i)\right\}>0
\end{equation}
for every $\beta\in(0,\infty)^d$.
\end{corollary}

\begin{proof}
By \eqref{eq:smoothed-defect},
$Q_{ab}n_c=n_c-h_c$ with $0<h_c<n_c\le M$.
The pointwise comparison
\eqref{eq:pointwise-risk-improvement} therefore gives
\[
 A(u)\ge-\log\left(1-\frac{h_c(u)}{n_c(u)}\right)
       \ge\frac{h_c(u)}{n_c(u)}
       \ge\frac{h_c(u)}{M},
\]
where we used $-\log(1-s)\ge s$ for $0<s<1$.

Let $L_c$ and $L_c'$ be independent and identically distributed.
Then
\begin{align*}
 \Delta(\beta)
 &\ge M^{-1}\E h_c(\theta+L_c)
 =M^{-1}\E e^{-\norm{\theta-\theta_0+L_c-L_c'}^2}.
\end{align*}
Since $L_c-L_c'$ has mean zero and
$\E\norm{L_c-L_c'}^2=2V_c$, we have
\[
 \E\norm{\theta-\theta_0+L_c-L_c'}^2
 =\norm{\theta-\theta_0}^2+2V_c.
\]
Jensen's inequality for the convex function $s\mapsto e^{-s}$
now gives
\[
 \Delta(\beta)\ge
 M^{-1}\exp\{-\norm{\theta-\theta_0}^2-2V_c\},
\]
which proves \eqref{eq:quantitative-bound}.
\end{proof}

\section{A three-dimensional exponential example}\label{sec:example}

Consider $d=3$ and $a_i=b_i=1$ for $i=1,2,3$.
The observation and future variables in coordinate $i$
are independent exponential variables with scale $\beta_i$;
the three scales remain unrestricted.
Since $B_i$ is uniform on $(0,1)$, $D_i=-\log B_i$
is exponential with mean and variance one.
The difference $W_i=D_i'-D_i$ has density
\[
 k_1(z)=\tfrac12e^{-|z|}.
\]
Indeed, for $z\ge0$,
\[
\int_0^\infty e^{-v}e^{-(v+z)}\dd v
=\frac12e^{-z},
\]
and the density is symmetric.
Consequently,
\begin{equation*}\label{eq:example-kernel}
 Q_{ab}h(u)=\frac18\int_{\R^3}h(u+z)e^{-\norm{z}_1}\dd z,
 \qquad
 \phi_{ab}(t)=\prod_{i=1}^3(1+t_i^2)^{-1},
\end{equation*}
where $\norm{z}_1=|z_1|+|z_2|+|z_3|$.
In this case $K=I_3$.
Proposition~\ref{prop:potential-tail} therefore gives
\[
 w_{ab}(\theta_0+v)\sim\frac{\sqrt\pi}{4\norm{v}}
 \qquad\text{as }\norm{v}\longrightarrow\infty.
\]

\subsection{The prior}

For any fixed $\theta_0\in\R^3$, the potential has the direct series
representation
\begin{align*}
 w_{ab}(u)={}&e^{-\norm{u-\theta_0}^2}
 +\sum_{n=1}^{\infty}8^{-n}
    \int_{(\R^3)^n}
       e^{-\norm{u+\sum_{j=1}^nz_j-\theta_0}^2
                    -\sum_{j=1}^n\norm{z_j}_1}
       \dd z_1\cdots\dd z_n.
       \label{eq:example-potential}
\end{align*}
All terms in the series are nonnegative. Alternatively,
\begin{equation}\label{eq:example-Fourier}
 w_{ab}(u)=\frac{\pi^{3/2}}{(2\pi)^3}
 \int_{\R^3}
 \frac{\cos\{t\cdot(u-\theta_0)\}e^{-\norm{t}^2/4}}
      {1-\prod_{i=1}^3(1+t_i^2)^{-1}}\dd t.
\end{equation}
Absolute convergence follows from the general result
\eqref{eq:potential-fourier}.

The corresponding prior is
\begin{equation}\label{eq:example-prior}
 \pi_{w_{ab}}(\beta)\dd\beta
 =\frac{w_{ab}(\log\beta_1,\log\beta_2,\log\beta_3)}
        {\beta_1\beta_2\beta_3}\dd\beta_1\dd\beta_2\dd\beta_3.
\end{equation}
It has infinite total mass, while each posterior is proper.

\subsection{Predictive densities and risk improvement}

Write $n_1$ for $n_{(1,1,1)}$ and $n_2$ for $n_{(2,2,2)}$.
The two smoothing integrals are
\begin{align*}
 n_1(u)&=\int_{(0,\infty)^3}
        w_{ab}(u-\log g)e^{-\sum_i g_i}\dd g, \\ 
 n_2(u)&=\int_{(0,\infty)^3}
        w_{ab}(u-\log g)\prod_{i=1}^3(g_i e^{-g_i})\dd g.
        \label{eq:example-n2}
\end{align*}
The relation $n_a=P_{ab+}n_c$ takes the explicit form
\[
 n_1(u)=\int_{(0,\infty)^3}n_2(u+v)e^{-\sum_i v_i}\dd v.
\]
The predictive densities are
\begin{equation}\label{eq:example-predictive}
 q_J(y\mid x)=\prod_{i=1}^3\frac{x_i}{(x_i+y_i)^2},
 \qquad
 q_{w_{ab}}(y\mid x)=q_J(y\mid x)
              \frac{n_2(\log(x+y))}{n_1(\log x)}.
\end{equation}
Both are normalized for every positive $x$.

The risk improvement can also be expressed as an absolutely
convergent iterated integral.
The density of $L_2=L_{(2,2,2)}$ is
\[
 \ell_2(v)=\exp\left\{2\sum_{i=1}^3v_i-\sum_{i=1}^3e^{v_i}\right\},
 \qquad v\in\R^3,
\]
as follows by $g_i=e^{v_i}$ in the Gamma density. Hence
\begin{align*}
 \Delta(\beta)
 =\int_{\R^3}\ell_2(v)\Bigg[&\log n_2(\theta+v)
 -\int_{(0,\infty)^3}e^{-\sum_i z_i}
             \log n_1(\theta+v-z)\dd z\Bigg]\dd v.
             \label{eq:example-risk-integral}
\end{align*}
The absolute convergence of both logarithmic terms follows from
Lemma~\ref{lem:log-bounds} and the exact coupling, as in
\eqref{eq:A-integrable}.

Substituting $a_i=b_i=1$ in \eqref{eq:CJ} gives
\[
 R(\beta,q_J)=3(1-\gamma),
 \qquad
 R(\beta,q_{w_{ab}})=R(\beta,q_J)-\Delta(\beta).
\]
Here $\gamma=-\psi(1)$ is Euler's constant and
$\psi(2)=\psi(1)+1$.

\subsection{Numerical illustration along a common-scale slice}
\label{sec:numerical}

For the numerical illustration we set $\theta_0=0$, so the shrinkage
target is $\beta_0=(1,1,1)$, and evaluate risks at
$\beta_1=\beta_2=\beta_3=\beta$. The predictive densities remain the
three-parameter rules \eqref{eq:example-predictive}; the equality of the
three true scales specifies the points at which risk is evaluated.
The same prior is used throughout the slice.

\begin{figure}[t]
\centering
\includegraphics[width=0.95\linewidth]{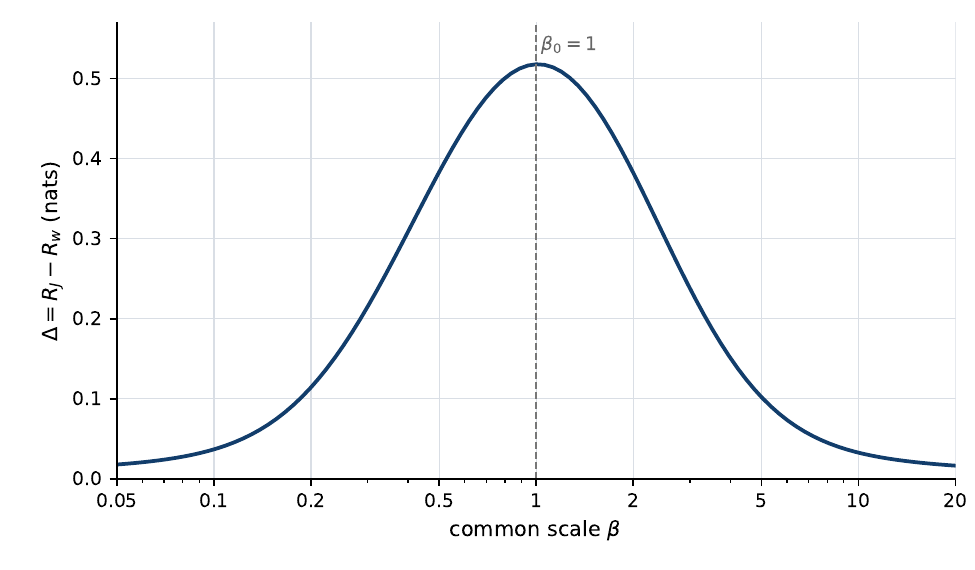}
\caption{Computed Kullback--Leibler risk improvement for $d=3$,
$a_i=b_i=1$, and $\theta_0=0$, evaluated at
$\beta_1=\beta_2=\beta_3=\beta$.
The vertical marker identifies the target $\beta=1$.
The curve shows the risk improvement of the predictive density
based on the prior in \eqref{eq:example-prior}.}
\label{fig:risk-improvement}
\end{figure}

Figure~\ref{fig:risk-improvement} shows the computed improvement
$\Delta(\beta)=R(\beta,q_J)-R(\beta,q_{w_{ab}})$ for $0.05\le\beta\le20$,
with a logarithmic horizontal axis. Table~\ref{tab:risks} gives selected
values. At the target $\beta=1$, the computed risk is approximately
$0.751$, compared with the Jeffreys risk
$3(1-\gamma)\approx1.268$, a reduction of about $40.8\%$.
The largest computed improvement on the 107-point evaluation grid
occurs at $\beta=1$; the improvement is smaller toward both ends of
the displayed interval.

\begin{table}[t]
\caption{Selected computed risks in the exponential example, in nats.
The Jeffreys risk is approximately $1.268$ throughout.
The percentage reduction is $100\Delta/R(\beta,q_J)$.}
\label{tab:risks}
\centering
\begin{tabular}{rrrr}
\hline
$\beta$ & $\Delta(\beta)$ & $R(\beta,q_{w_{ab}})$ & Reduction (\%)\\
\hline
$0.050$ & $0.018$ & $1.251$ & $1.406$\\
$0.200$ & $0.114$ & $1.154$ & $9.009$\\
$1.000$ & $0.518$ & $0.751$ & $40.828$\\
$5.000$ & $0.102$ & $1.166$ & $8.047$\\
$20.000$ & $0.016$ & $1.252$ & $1.292$\\
\hline
\end{tabular}
\end{table}

These values were obtained by deterministic Fourier quadrature for
the potential and smoothed logarithms. Domain and spacing refinements,
together with a separate scrambled-Sobol calculation of the outer
expectation, are described in Appendix~\ref{app:numerical}.
At the reciprocal scales $\beta=0.2$ and $\beta=5$, the computed
improvements differ by $0.012$. Evenness of the prior weight does not
imply reciprocal-scale risk symmetry: the log-Gamma sampling law is
asymmetric, and the risk averages the logarithm of a smoothed weight.
This difference persists under the numerical checks in that appendix.
The illustration quantifies improvement on this slice; strict domination
on the full parameter space follows from Theorem~\ref{thm:main}.

\section{Discussion}\label{sec:discussion}

The constructed shrinkage predictive densities strictly dominate the
Jeffreys rule and are minimax for arbitrary known positive shape vectors
when $d\ge3$. Their target $\beta_0$ can be specified in advance without
restricting the unknown coordinate scales.

The beta--gamma coupling identifies a discrete transition determined
by the observation and future distributions.
The predictive comparison uses $A_a=P_{ab+}A_c$ and
$Q_{ab}=P_{ab-}P_{ab+}$, while the potential satisfies
$(I-Q_{ab})w_{ab}=f$.
These identities explain why the full distribution of the increment
is relevant to the sufficient condition.
The covariance of the increment determines the quadratic behavior
of its characteristic function near the origin, while the exact
risk comparison uses the full transition distribution.

The prior depends on the known observation and future shapes
through $Q_{ab}$.
For each such design, the resulting predictive density
improves $q_J$ over the entire unrestricted scale parameter space.
The center is a fixed location in logarithmic scale coordinates,
and changing it translates the potential without changing $M$.
The tail equivalent in Proposition~\ref{prop:potential-tail}
connects this construction to the classical harmonic power $2-d$,
expressed in the quadratic form determined by $K^{-1}$.
The bounds establish the construction for arbitrarily small positive
shapes, although their constants need not remain uniform when the
shapes approach the boundary of the positive orthant.

The Gaussian form of the source is not essential for strict
domination or minimaxity. Its use here gives a bounded continuous
weight that is symmetric about a prespecified center, has its
unique maximum there, and admits the explicit tail equivalent
in Proposition~\ref{prop:potential-tail}. These additional
properties need not hold for a general source.

More generally, let $\nu$ be a nonnegative measure on $\R^d$
with $0<\nu(\R^d)<\infty$, and define the potential measure
\[
 w_{ab}^{\nu}
 =\sum_{n=0}^{\infty}Q_{ab}^n\nu,
\]
where $Q_{ab}$ acts on measures by convolution with the
distribution of $D'_{ab}-D_{ab}$.
This gives a prior measure in logarithmic scale coordinates.
Let $\ell_r$ denote the density of $L_r$, and put
\[
 n_r^{\nu}(u)
 =\int_{\R^d}\ell_r(u-t)w_{ab}^{\nu}(\dd t),
 \qquad r\in\{a,c\}.
\]
For $d\ge3$, these functions are positive and bounded and
satisfy the required logarithmic integrability conditions.
Moreover,
\[
 n_c^{\nu}(u)-Q_{ab}n_c^{\nu}(u)
 =\int_{\R^d}\ell_c(u-t)\nu(\dd t)>0.
\]
The predictive density $q_{w_{ab}^{\nu}}$, obtained from
\eqref{eq:qw} by replacing $n_r$ with $n_r^{\nu}$, therefore
strictly dominates $q_J$ and is minimax.
This includes sources $\nu(\dd t)=f(t)\dd t$ with
$f\ge0$ and $0<\int_{\R^d}f(t)\dd t<\infty$,
as well as the Dirac source $\nu=\delta_{\theta_0}$.
A full justification of this extension, including posterior
normalization and logarithmic integrability, will be given
elsewhere.

Multiplying the source by a positive constant does not change
the resulting predictive density. In particular, the Gaussian
source may be normalized to a probability density and may
have any positive variance. As its variance tends to zero,
the normalized source converges weakly to
$\delta_{\theta_0}$. The Dirac-source result is established
directly through the smoothed potential identity above,
rather than by passing the strict risk inequality to this
limit.

Beyond Gamma models, the Markov potential construction may apply
to other prediction problems in which an exact coupling relates
the observed and combined-data experiments through an
independent random translation. Such a coupling can identify
a transition operator whose superharmonic functions yield
predictive risk improvements. Normal location models provide
another example of this structure.

\begin{acks}[Acknowledgments]
Large language model (LLM) tools were used in the research and
writing of this paper. The Markov potential prior presented here was
proposed and derived with the assistance of these tools.
The authors take full responsibility for the content of the manuscript.

This work was supported by JSPS KAKENHI Grant Numbers 22H00510 and 25H01156.
\end{acks}

\begin{appendix}

\section{The weighted posterior and predictive density}
\label{app:distributional}

\begin{proof}[Proof of Lemma~\ref{lem:posterior}]
Choose $C<\infty$ such that $w\le C$.
In one coordinate, the change of variables $g=x_i/\beta_i$ has
absolute Jacobian $x_i g^{-2}$ and gives
\[
 \frac{x_i^{a_i-1}}{\Gamma(a_i)}
       \beta_i^{-a_i-1}e^{-x_i/\beta_i}\dd\beta_i
 =x_i^{-1}\frac{g^{a_i-1}e^{-g}}{\Gamma(a_i)}\dd g.
\]
Consequently, the Jeffreys marginal and posterior are
\begin{align}
 m_J(x)&=\prod_i x_i^{-1},
 \label{eq:Jeffreys-marginal}\\
 \pi_J(\beta\mid x)
 &=\prod_i\frac{x_i^{a_i}}{\Gamma(a_i)}
      \beta_i^{-a_i-1}e^{-x_i/\beta_i}.
 \label{eq:Jeffreys-posterior}
\end{align}

Applying the same change of variables to the weighted prior
in \eqref{eq:prior} gives
\begin{align}
 m_w(x)
 &=m_J(x)\int_{(0,\infty)^d}
 w(\log x-\log g)
 \prod_i\frac{g_i^{a_i-1}e^{-g_i}}{\Gamma(a_i)}
 \,\dd g \notag\\
 &=m_J(x)n_a(\log x).
 \label{eq:weighted-marginal}
\end{align}
By Bayes' formula,
\[
 \pi_w(\beta\mid x)
 =\frac{w(\log\beta)}{n_a(\log x)}
 \pi_J(\beta\mid x).
\]
Since $0<n_a(u)\le C$ for every $u$, the posterior is proper
for every positive $x$.

For the joint marginal, the substitution
$g_i=(x_i+y_i)/\beta_i$ gives
\begin{align}
 m_J(x,y)
 &=\prod_i\frac{\Gamma(c_i)}{\Gamma(a_i)\Gamma(b_i)}
      \frac{x_i^{a_i-1}y_i^{b_i-1}}{(x_i+y_i)^{c_i}},
       \label{eq:Jeffreys-joint}\\
 m_w(x,y)
 &=m_J(x,y)n_c(\log(x+y)).
 \label{eq:weighted-joint}
\end{align}
Taking the ratio $m_w(x,y)/m_w(x)$ gives
the predictive density \eqref{eq:qw}.

Finally, the predictive density $q_w$ is normalized
because it is a mixture of the densities $p_b(y\mid\beta)$
with respect to the proper posterior
$\pi_w(\beta\mid x)$.
The same argument applies to $q_J$.
\end{proof}

\section{Increment distribution and frequency bounds}
\label{app:kernel}

This appendix proves the properties of the increment used in
Section~\ref{sec:potential}. The arguments apply to every positive
integer dimension and all positive shapes.

\begin{lemma}\label{lem:kernel}
Let $B_{a_i b_i}\sim\operatorname{Beta}(a_i,b_i)$ be independent,
and define
$D_{ab}=(D_{a_i b_i})_{i=1}^d$ with
$D_{a_i b_i}=-\log B_{a_i b_i}$.
Let $D_{ab}'$ be an independent copy of $D_{ab}$,
and set $W_{ab}=D_{ab}'-D_{ab}$.
Then $\E\norm{W_{ab}}^2<\infty$,
$W_{ab}$ and $-W_{ab}$ have the same law,
$\E W_{ab}=0$, and
$\Cov(W_{ab})=2K$, where $K$ is the positive definite
matrix in \eqref{eq:K}.

The characteristic function $\phi_{ab}$ defined in
\eqref{eq:characteristic} satisfies
\begin{equation}\label{eq:phi-properties}
 0\le\phi_{ab}(t)\le1,\qquad
 \phi_{ab}(t)<1\ (t\ne0),\qquad
 \phi_{ab}(t)\longrightarrow0\
 (\norm{t}\longrightarrow\infty).
\end{equation}
Moreover,
\begin{equation}\label{eq:phi-expansion}
 1-\phi_{ab}(t)
 =t^{\mathsf T}Kt+o(\norm{t}^2),
 \qquad t\to0.
\end{equation}
The constants in \eqref{eq:frequency-bounds} exist
and can be chosen so that
$1-\phi_{ab}(t)\ge\kappa\norm{t}^2$
for $\norm{t}\le\delta$.
\end{lemma}

\begin{proof}
Changing variables $B_{a_i b_i}=e^{-v}$ in the Beta density gives
\begin{equation}\label{eq:D-density}
 g_i(v)=
 \frac{e^{-a_iv}(1-e^{-v})^{b_i-1}}
      {\mathrm B(a_i,b_i)}
 \one_{(0,\infty)}(v).
\end{equation}
The moment generating function of $D_{a_i b_i}$ is
\begin{equation}\label{eq:D-mgf}
 \E e^{sD_{a_i b_i}}
 =\frac{\mathrm B(a_i-s,b_i)}{\mathrm B(a_i,b_i)},
 \qquad s<a_i.
\end{equation}
Since this is finite in a neighborhood of zero,
all moments of $D_{a_i b_i}$ are finite.
Differentiating the logarithm of this moment generating
function at zero gives
\begin{equation}\label{eq:D-moments}
 \E D_{a_i b_i}=\psi(c_i)-\psi(a_i),
 \qquad
 \Var(D_{a_i b_i})
 =\psi_1(a_i)-\psi_1(c_i)=K_i>0.
\end{equation}
Independence of coordinates and copies gives
$\E W_{ab}=0$ and $\Cov(W_{ab})=2K$.
Interchanging the two copies proves symmetry.

The characteristic function of $D_{a_i b_i}$ is
\begin{equation}\label{eq:D-characteristic}
 \chi_{ab,i}(s):=\E e^{isD_{a_i b_i}}
 =\frac{\mathrm B(a_i-is,b_i)}{\mathrm B(a_i,b_i)}
 =\frac{\Gamma(a_i-is)\Gamma(c_i)}
        {\Gamma(a_i)\Gamma(c_i-is)}.
\end{equation}
By independence,
\[
 \phi_{ab}(t)=\prod_i|\chi_{ab,i}(t_i)|^2\in[0,1].
\]
Since $D_{a_i b_i}$ has a density,
$|\chi_{ab,i}(s)|<1$ for $s\ne0$.
Hence $\phi_{ab}(t)<1$ for $t\ne0$.
Moreover, the Riemann--Lebesgue lemma gives
$\chi_{ab,i}(s)\to0$ as $|s|\to\infty$.
Since there are finitely many coordinates,
$\phi_{ab}(t)\to0$ as $\norm{t}\to\infty$.

By symmetry and the finite second moment of $W_{ab}$,
\[
 1-\phi_{ab}(t)
 =\tfrac12\E(t\cdot W_{ab})^2+o(\norm{t}^2)
 =t^{\mathsf T}Kt+o(\norm{t}^2),
\]
which proves \eqref{eq:phi-expansion}.
Since $K$ is positive definite, there exist
$\delta,\kappa>0$ such that
$1-\phi_{ab}(t)\ge\kappa\norm{t}^2$
for $\norm{t}\le\delta$.
Consequently,
\begin{equation}\label{eq:local-phi-bound}
 0\le\phi_{ab}(t)
 \le1-\kappa\norm{t}^2
 \le e^{-\kappa\norm{t}^2},
 \qquad \norm{t}\le\delta.
\end{equation}

Finally, since $\phi_{ab}(t)<1$ for $t\ne0$
and $\phi_{ab}(t)\to0$ as $\norm{t}\to\infty$,
continuity gives a constant $0<\rho<1$ such that
$\phi_{ab}(t)\le\rho$ whenever $\norm{t}\ge\delta$.
This proves \eqref{eq:frequency-bounds}.
\end{proof}

\section{Pointwise Fourier representation and the potential integral}
\label{app:Fourier}

We use the Fourier convention
\[
 \widehat h(t)=\int_{\R^d}e^{-it\cdot u}h(u)\dd u.
\]

\begin{proof}[Proof of Lemma~\ref{lem:iterate-fourier}]
The Gaussian Fourier inversion formula gives
\begin{equation}\label{eq:gaussian-inverse}
 f(u)=(2\pi)^{-d}\int_{\R^d}
 e^{it\cdot(u-\theta_0)}
 \pi^{d/2}e^{-\norm{t}^2/4}\dd t.
\end{equation}
Let $W_{ab,1},W_{ab,2},\ldots$ be independent copies
of $W_{ab}$, and put
$S_n=\sum_{j=1}^nW_{ab,j}$, with $S_0=0$.
Then
\[
 Q_{ab}^nf(u)=\E f(u+S_n),
 \qquad
 \E e^{it\cdot S_n}=\phi_{ab}(t)^n.
\]
Since the Gaussian Fourier transform is integrable,
Fubini's theorem allows us to take expectations in
\eqref{eq:gaussian-inverse}, giving
\begin{align*}
 Q_{ab}^nf(u)
 &=(2\pi)^{-d}\int_{\R^d}
 e^{it\cdot(u-\theta_0)}
 \pi^{d/2}e^{-\norm{t}^2/4}
 \E e^{it\cdot S_n}\dd t\\
 &=(2\pi)^{-d}\int_{\R^d}
 e^{it\cdot(u-\theta_0)}
 \pi^{d/2}e^{-\norm{t}^2/4}
 \phi_{ab}(t)^n\dd t.
\end{align*}
This proves \eqref{eq:Qn-pointwise-fourier}.
Absolute convergence follows from
$|\phi_{ab}(t)|\le1$.
\end{proof}

To prove \eqref{eq:potential-fourier}, suppose $d\ge3$ and put
\[
 F_{ab}(t)=
 \frac{\pi^{d/2}e^{-\norm{t}^2/4}}
      {1-\phi_{ab}(t)},\qquad t\ne0.
\]
By Lemma~\ref{lem:kernel}, for
$0<\norm{t}\le\delta$,
\[
 F_{ab}(t)\le
 \frac{\pi^{d/2}}{\kappa\norm{t}^2}.
\]
This bound is integrable near zero when $d>2$,
since its radial integral is proportional to
$\int_0^\delta r^{d-3}\dd r$.
For $\norm{t}\ge\delta$, we have
\[
 F_{ab}(t)\le
 \frac{\pi^{d/2}e^{-\norm{t}^2/4}}{1-\rho},
\]
which is integrable. Thus $F_{ab}$ is integrable
over $\R^d$, regardless of its value at the origin.

By Lemma~\ref{lem:iterate-fourier}, the finite sums
$\sum_{n=0}^N Q_{ab}^nf$ have Fourier integral
representations with multipliers
$\sum_{n=0}^N\phi_{ab}(t)^n$.
Since
\[
 0\le\sum_{n=0}^N\phi_{ab}(t)^n
 \le\frac{1}{1-\phi_{ab}(t)}
 \qquad(t\ne0),
\]
the integrands are bounded in absolute value by the integrable
function $F_{ab}$.
Dominated convergence therefore allows us to let $N\to\infty$
under the integral.
Using
\[
 \sum_{n=0}^\infty\phi_{ab}(t)^n
 =\frac{1}{1-\phi_{ab}(t)}
 \qquad(t\ne0),
\]
and the convergence of the potential series,
we obtain
\[
 w_{ab}(u)=(2\pi)^{-d}\int_{\R^d}
 e^{it\cdot(u-\theta_0)}F_{ab}(t)\dd t.
\]
Finally, $F_{ab}$ is even, so the sine term
integrates to zero. This gives the cosine
representation \eqref{eq:potential-fourier}.

\subsection{The tail equivalent}\label{app:tail}

We prove Proposition~\ref{prop:potential-tail}.
Throughout the proof, the shapes and hence the positive definite
matrix $K$ are fixed. Put $v=u-\theta_0$ and
$k_K(t)=t^{\mathsf T}Kt$. For $t\ne0$, define
\[
 F_{0,ab}(t)=
 \frac{\pi^{d/2}e^{-\norm{t}^2/4}}{k_K(t)},
 \qquad
 \mathcal R_{ab,K}(t)=F_{ab}(t)-F_{0,ab}(t),
\]
where $F_{ab}$ is defined above.
Both $F_{ab}$ and $F_{0,ab}$ are integrable when $d>2$.
Writing
\[
 H_{ab,K}(v)=(2\pi)^{-d}\int_{\R^d}
       e^{it\cdot v}F_{0,ab}(t)\dd t,
\]
and
\[
 \mathcal E_{ab,K}(v)=(2\pi)^{-d}\int_{\R^d}
       e^{it\cdot v}\mathcal R_{ab,K}(t)\dd t,
\]
we have
\[
 w_{ab}(\theta_0+v)=H_{ab,K}(v)+\mathcal E_{ab,K}(v).
\]

\subsubsection{The Fourier remainder}

By Appendix~\ref{app:kernel}, $W_{ab}$ has moments of
all orders. Hence $\phi_{ab}$ is infinitely differentiable,
with bounded derivatives of every order.
Moreover, $\phi_{ab}$ is even and $\Cov(W_{ab})=2K$.
Taylor's formula therefore gives
\[
 1-\phi_{ab}(t)=k_K(t)+e_{ab}(t),
 \qquad
 |\partial^\alpha e_{ab}(t)|
 \le C_\alpha\norm{t}^{4-|\alpha|}
\]
for $0<\norm{t}\le\delta_0$, where $\delta_0>0$ is
sufficiently small and $\alpha$ is any fixed multiindex.

Since both $k_K(t)$ and $1-\phi_{ab}(t)$ are bounded
above and below by positive constant multiples of
$\norm{t}^2$ near zero, differentiation of
\[
 \mathcal R_{ab,K}(t)
 =-\pi^{d/2}e^{-\norm{t}^2/4}
   \frac{e_{ab}(t)}
        {k_K(t)(1-\phi_{ab}(t))}
\]
gives
\begin{equation}\label{eq:tail-remainder-derivatives}
 |\partial^\alpha\mathcal R_{ab,K}(t)|
 \le C_\alpha\norm{t}^{-|\alpha|}
 \qquad(0<\norm{t}\le\delta_0).
\end{equation}
Away from zero, all derivatives of
$\mathcal R_{ab,K}$ are integrable because
$1-\phi_{ab}$ is bounded away from zero and the
Gaussian factor controls the tails.
Consequently, the derivatives of orders up to $d-1$
are integrable.

These derivatives are also weak derivatives on $\R^d$.
Indeed, integration by parts outside a ball of radius
$\varepsilon$ produces boundary terms of order
$O(\varepsilon^{d-j})$ at derivative order $j$,
which vanish for $j\le d-1$.

Fourier integration by parts thus gives, for
$j=1,\ldots,d$,
\[
 |v_j|^{d-1}|\mathcal E_{ab,K}(v)|
 \le(2\pi)^{-d}
 \norm{\partial_j^{d-1}\mathcal R_{ab,K}}_{L^1}.
\]
Since some coordinate satisfies
$|v_j|\ge\norm{v}/\sqrt d$, we obtain
\begin{equation}\label{eq:tail-remainder-decay}
 \mathcal E_{ab,K}(v)
 =O(\norm{v}^{1-d})
 =o(\norm{v}^{2-d}),
\end{equation}
uniformly in direction.

\subsubsection{The leading inverse transform}

Using
\[
 \frac{1}{k_K(t)}
 =\int_0^\infty e^{-s k_K(t)}\dd s
\]
and the Gaussian Fourier formula gives
\begin{equation}\label{eq:tail-heat-representation}
 H_{ab,K}(v)=\int_0^\infty
 \frac{\exp\{-v^{\mathsf T}(I+4sK)^{-1}v\}}
      {\sqrt{\det(I+4sK)}}\dd s.
\end{equation}
The interchange of integrals is justified by
\[
 \int_{\R^d}
 \frac{e^{-\norm{t}^2/4}}{k_K(t)}\dd t<\infty,
\]
which holds for $d>2$.

Write $v=rz$, where $r=\norm{v}\ge1$ and
$\norm{z}=1$. Substituting $s=r^2\tau$ gives
\[
 r^{d-2}H_{ab,K}(rz)
 =\int_0^\infty J_r(\tau,z)\dd\tau,
\]
where
\[
 J_r(\tau,z)=
 \frac{\exp\{-z^{\mathsf T}
              (r^{-2}I+4\tau K)^{-1}z\}}
      {\sqrt{\det(r^{-2}I+4\tau K)}}.
\]
For each $\tau>0$, $J_r(\tau,z)$ converges
uniformly over $\norm{z}=1$ to
\[
 J_\infty(\tau,z)=
 \frac{(4\tau)^{-d/2}}{\sqrt{\det K}}
 \exp\left\{
 -\frac{z^{\mathsf T}K^{-1}z}{4\tau}
 \right\}.
\]

To justify integration of the limit, let
$c_K=(1+4\lambda_{\max}(K))^{-1}>0$.
For $0<\tau<r^{-2}$,
\[
 J_r(\tau,z)\le r^d e^{-c_Kr^2},
\]
so the integral over this interval tends to zero
uniformly in $z$.
For $\tau\ge r^{-2}$,
\[
 J_r(\tau,z)\le C\tau^{-d/2}e^{-c_K/\tau},
\]
which is integrable on $(0,\infty)$ when $d>2$.
The same bound holds for $J_\infty(\tau,z)$.
For each $\tau>0$,
\[
 \sup_{\norm{z}=1}
 \left|
 \one_{\{\tau\ge r^{-2}\}}J_r(\tau,z)
 -J_\infty(\tau,z)
 \right|\longrightarrow0.
\]
This expression is bounded by
$2C\tau^{-d/2}e^{-c_K/\tau}$, which is integrable.
Dominated convergence, together with the preceding estimate
on $(0,r^{-2})$, therefore gives
\[
 r^{d-2}H_{ab,K}(rz)
 \longrightarrow
 \int_0^\infty J_\infty(\tau,z)\dd\tau
\]
uniformly over $\norm{z}=1$.

Writing $\eta=z^{\mathsf T}K^{-1}z$ and substituting
$\zeta=\eta/(4\tau)$, we obtain
\[
\begin{aligned}
 \int_0^\infty J_\infty(\tau,z)\dd\tau
 &=\frac{\eta^{1-d/2}}{4\sqrt{\det K}}
   \int_0^\infty
   \zeta^{d/2-2}e^{-\zeta}\dd\zeta\\
 &=\frac{\Gamma(d/2-1)}{4\sqrt{\det K}}
   (z^{\mathsf T}K^{-1}z)^{(2-d)/2}.
\end{aligned}
\]
This limit is continuous and bounded away from zero
on the unit sphere. Combining the uniform limit with
\eqref{eq:tail-remainder-decay} proves
\eqref{eq:potential-tail}, uniformly in direction.
Finally, positive definiteness of $K$ implies that
$v^{\mathsf T}K^{-1}v$ is comparable to $\norm{v}^2$,
giving the claimed norm order.

\section{Logarithmic integrability and risk calculations}
\label{app:integrability}

We first establish the logarithmic moment bounds needed for
the risk comparison. We then show that the logarithmic losses
for both predictive densities are absolutely integrable,
justifying the subtraction of their risks.

For any probability densities $p$ and $q$, the inequality
$z\log(1/z)\le e^{-1}$ for $0<z<1$ gives
\begin{equation}\label{eq:KL-negative-part}
 \int p(y)\left[\log\frac{p(y)}{q(y)}\right]^-\dd y
 \le e^{-1}\int q(y)\dd y=e^{-1}.
\end{equation}
Here $s^-=\max\{-s,0\}$, with the usual extended
conventions at zeros.
Applied conditionally on $X$, this bound ensures that
\eqref{eq:risk} is well defined for every $q\in\mathcal Q$,
even when the risk is infinite.

\subsection{Log-Gamma moments and the smoothing envelope}

For $G\sim\operatorname{Gamma}(r,1)$ with $r>0$,
all absolute logarithmic moments are finite:
\begin{equation}\label{eq:logGamma-integral}
 \E|\log G|^k
 =\frac1{\Gamma(r)}\int_0^\infty
       g^{r-1}e^{-g}|\log g|^k\dd g<\infty
\end{equation}
for every nonnegative integer $k$.
Indeed, integrability near zero follows by substituting
$v=-\log g$, while the Gamma tail ensures integrability
at infinity.

Moreover,
\[
 \E G^s=\frac{\Gamma(r+s)}{\Gamma(r)},
 \qquad s>-r.
\]
Since the moment generating function of $\log G$ is finite
in a neighborhood of zero, differentiating its logarithm
at zero gives
\[
 \E\log G=\psi(r),\qquad
 \Var(\log G)=\psi_1(r).
\]
Summing over the independent coordinates proves
\eqref{eq:log-moments}.

\begin{proof}[Proof of Lemma~\ref{lem:log-bounds}]
The preceding calculation establishes the moments.
Since $w_{ab}\ge f$, Jensen's inequality gives
\begin{align*}
 n_r(u)
 &\ge\E e^{-\norm{u-\theta_0-L_r}^2}\\
 &\ge\exp\{-\E\norm{u-\theta_0-L_r}^2\}\\
 &=\exp\{-\norm{u-\theta_0-\mu_r}^2-V_r\}.
\end{align*}
This proves \eqref{eq:smoothing-lower}.
Combining this lower bound with $n_r\le M$ and
$M\ge1$ proves \eqref{eq:log-envelope}.

Substituting $u=\theta+L_r$ gives
\begin{align*}
 \E|\log n_r(\theta+L_r)|
 &\le\log M+
    \E\norm{\theta-\theta_0+L_r-\mu_r}^2+V_r\\
 &=\log M+\norm{\theta-\theta_0}^2+2V_r,
\end{align*}
which proves \eqref{eq:outer-log-bound}.

Applying \eqref{eq:log-envelope} to
$n_a(u-D_{ab})$ proves \eqref{eq:inner-log-bound}.
Indeed, since $\E D_{ab}=\mu_c-\mu_a$,
\[
 \E\norm{u-D_{ab}-\theta_0-\mu_a}^2
 =\norm{u-\theta_0-\mu_c}^2
  +\sum_{i=1}^d K_i<\infty.
\]
\end{proof}

\subsection{Absolute integrability of the logarithmic losses}

\begin{lemma}\label{lem:original-losses}
For every positive $a,b$ and every
$\beta\in(0,\infty)^d$, the logarithmic loss
$\log\{p_b(Y\mid\beta)/q_J(Y\mid X)\}$
is absolutely integrable. Its risk is the finite constant
\begin{equation}\label{eq:Jeffreys-risk}
 R(\beta,q_J)=\sum_{i=1}^d
 \bigl[\log\Gamma(a_i)-\log\Gamma(c_i)-b_i
       -a_i\psi(a_i)+c_i\psi(c_i)\bigr].
\end{equation}
If $w$ is positive, bounded, and Borel measurable
and satisfies \eqref{eq:criterion-outer}, then
the logarithmic loss for $q_w$ is also absolutely
integrable. Both risks are nonnegative.
\end{lemma}

\begin{proof}
Write $U_i=X_i/\beta_i$ and $V_i=Y_i/\beta_i$.
These are independent unit-scale Gamma variables with
shapes $a_i$ and $b_i$, and $U_i+V_i$ has shape $c_i$.
Direct calculation gives
\begin{equation}\label{eq:original-Jeffreys-loss}
 \log\frac{p_b(Y\mid\beta)}{q_J(Y\mid X)}
 =\sum_{i=1}^d
 \bigl[\log\Gamma(a_i)-\log\Gamma(c_i)-V_i
       -a_i\log U_i+c_i\log(U_i+V_i)\bigr].
\end{equation}
All terms are absolutely integrable by
\eqref{eq:logGamma-integral} and the finite first
moments of the Gamma variables.
Taking expectations gives \eqref{eq:Jeffreys-risk}.

Absolute integrability of the logarithmic loss for $q_w$
follows from the predictive-ratio bound established
below in Appendix~\ref{app:endpoint-regularity}.
Finally, both risks are nonnegative by Gibbs' inequality,
applied conditionally on $X$.
\end{proof}

Formula~\eqref{eq:Jeffreys-risk} gives the constant
risk used in the example and in the minimax calculation.

\subsection{Justification of the risk comparison}
\label{app:endpoint-regularity}

For a positive bounded weight $w$ satisfying
\eqref{eq:criterion-outer}, the predictive ratio
\eqref{eq:qw} gives
\[
 \E\left|\log\frac{q_w(Y\mid X)}{q_J(Y\mid X)}\right|
 \le\E|\log n_c(\theta+L_c)|
    +\E|\log n_a(\theta+L_a)|<\infty.
\]
This bound uses only the marginal distributions of
$X$ and $X+Y$, without requiring independence.

Since
\[
 \log\frac{p_b(Y\mid\beta)}{q_w(Y\mid X)}
 =\log\frac{p_b(Y\mid\beta)}{q_J(Y\mid X)}
  -\log\frac{q_w(Y\mid X)}{q_J(Y\mid X)},
\]
the absolute integrability of
\eqref{eq:original-Jeffreys-loss} also establishes
absolute integrability of the logarithmic loss for $q_w$.
This justifies the risk subtraction in
\eqref{eq:risk-marginals}.

For the potential $w_{ab}$,
Lemma~\ref{lem:log-bounds} gives
\begin{equation}\label{eq:explicit-ratio-log-bound}
 \E\left|
 \log\frac{q_{w_{ab}}(Y\mid X)}{q_J(Y\mid X)}
 \right|
 \le 2\log M+2\norm{\theta-\theta_0}^2+2V_a+2V_c.
\end{equation}

For the independent coupling of $L_c$ and $D_{ab}$
in \eqref{eq:coupling}, we have
\[
 \E_{L_c,D_{ab}}
 |\log n_a(\theta+L_c-D_{ab})|
 =\E|\log n_a(\theta+L_a)|<\infty.
\]
Thus Fubini's theorem gives \eqref{eq:endpoint-identity}.
If \eqref{eq:criterion-inner} also holds, the conditional
integral is finite for every $u\in\R^d$.
Moreover, the same coupling and the triangle inequality give
\begin{equation}\label{eq:A-integrable}
 \E|A(\theta+L_c)|
 \le\E|\log n_c(\theta+L_c)|
    +\E|\log n_a(\theta+L_a)|<\infty.
\end{equation}
Hence all expectations used in the risk comparison are finite.

\section{A proper-prior verification of the minimax lower bound}
\label{app:minimax}

We prove the minimax lower bound in Proposition~\ref{prop:minimax}
using a sequence of proper priors.
Fix $\tau_i>0$ and $\varepsilon>0$.
In rate coordinates $\lambda_i=1/\beta_i$, consider
the product prior
\begin{equation}\label{eq:minimax-prior-rate}
 \pi_\varepsilon(\dd\lambda)
 =\prod_{i=1}^d
   \frac{\tau_i^\varepsilon}{\Gamma(\varepsilon)}
   \lambda_i^{\varepsilon-1}e^{-\tau_i\lambda_i}
   \dd\lambda_i.
\end{equation}
Equivalently, its density in the scale coordinates is
\[
 \prod_{i=1}^d
 \frac{\tau_i^\varepsilon}{\Gamma(\varepsilon)}
 \beta_i^{-\varepsilon-1}e^{-\tau_i/\beta_i}.
\]
We use $\pi_\varepsilon$ for the same prior in either
parameterization.

\subsection{Predictive density and Bayes risk}

By Gamma conjugacy, the posterior rates are independent,
with shapes $a_i+\varepsilon$ and rates $x_i+\tau_i$.
Integrating the future Gamma density against this posterior gives
\begin{equation}\label{eq:epsilon-predictive}
 q_\varepsilon(y\mid x)=\prod_{i=1}^d
 \frac{\Gamma(c_i+\varepsilon)}
      {\Gamma(a_i+\varepsilon)\Gamma(b_i)}
 \frac{(x_i+\tau_i)^{a_i+\varepsilon}y_i^{b_i-1}}
      {(x_i+y_i+\tau_i)^{c_i+\varepsilon}}.
\end{equation}
As a posterior mixture, $q_\varepsilon(\cdot\mid x)$
is a proper predictive density for every positive $x$.
Let
\[
 m_\varepsilon(x)
 =\int p_a(x\mid\beta)\pi_\varepsilon(\dd\beta)
\]
be the corresponding sampling marginal.

Under the joint prior and sampling distribution, define
\[
 T_i=\tau_i\lambda_i,\qquad
 U_i=\lambda_iX_i,\qquad
 V_i=\lambda_iY_i.
\]
These variables are independent unit-scale Gamma variables
with respective shapes $\varepsilon,a_i,b_i$.
Consequently,
$T_i+U_i\sim\operatorname{Gamma}(a_i+\varepsilon,1)$
and
$T_i+U_i+V_i\sim\operatorname{Gamma}(c_i+\varepsilon,1)$.

Direct calculation gives
\begin{equation}\label{eq:epsilon-original-loss}
 \log\frac{p_b(Y\mid\beta)}{q_\varepsilon(Y\mid X)}
 =\sum_{i=1}^d\ell_{\varepsilon,i},
\end{equation}
where
\[
 \begin{aligned}
 \ell_{\varepsilon,i}
 ={}&\log\Gamma(a_i+\varepsilon)
     -\log\Gamma(c_i+\varepsilon)-V_i\\
 &-(a_i+\varepsilon)\log(T_i+U_i)
 +(c_i+\varepsilon)\log(T_i+U_i+V_i).
 \end{aligned}
\]
All terms are absolutely integrable by
\eqref{eq:logGamma-integral} and the finite first
moments of the Gamma variables.
Taking expectations gives the finite Bayes risk
\begin{equation}\label{eq:Bayes-risk-epsilon}
 \begin{aligned}
 r_\varepsilon=\sum_{i=1}^d\Bigl[
 &\log\Gamma(a_i+\varepsilon)
  -\log\Gamma(c_i+\varepsilon)-b_i\\
 &-(a_i+\varepsilon)\psi(a_i+\varepsilon)
  +(c_i+\varepsilon)\psi(c_i+\varepsilon)\Bigr].
 \end{aligned}
\end{equation}
By continuity of the Gamma and digamma functions,
$r_\varepsilon\to C_J$ as $\varepsilon\downarrow0$,
with the rates $\tau_i$ fixed.

\subsection{Comparison with arbitrary predictive densities}

Let $q\in\mathcal Q$, without any factorization assumption.
By \eqref{eq:KL-negative-part}, the negative part of
its logarithmic loss has conditional integral at most $e^{-1}$.
Thus its Bayes risk under the proper prior
$\pi_\varepsilon$ is well defined, possibly infinite.

The identity
\[
 \log\frac{p_b(Y\mid\beta)}{q(Y\mid X)}
 =
 \log\frac{p_b(Y\mid\beta)}{q_\varepsilon(Y\mid X)}
 +\log\frac{q_\varepsilon(Y\mid X)}{q(Y\mid X)}
\]
has an absolutely integrable first term by
\eqref{eq:epsilon-original-loss}.
Conditional on $X=x$, the predictive density under
$\pi_\varepsilon$ is $q_\varepsilon(\cdot\mid x)$.
The negative part of the second term is integrable because
\[
 \int_{\{q>q_\varepsilon\}}q_\varepsilon(y\mid x)
 \log\frac{q(y\mid x)}{q_\varepsilon(y\mid x)}
 \dd y\le e^{-1}.
\]
These bounds justify the decomposition
\begin{equation}\label{eq:full-class-Bayes-decomposition}
 \begin{aligned}
 \int R(\beta,q)\pi_\varepsilon(\dd\beta)
 ={}&r_\varepsilon
 +\int m_\varepsilon(x)\Biggl[
   \int q_\varepsilon(y\mid x)
   \log\frac{q_\varepsilon(y\mid x)}{q(y\mid x)}
   \dd y\Biggr]\dd x.
 \end{aligned}
\end{equation}
The inner integral is a KL divergence and is nonnegative,
possibly infinite. The decomposition remains valid when
the Bayes risk is infinite.

Consequently,
\[
 \sup_{\beta\in\mathcal B}R(\beta,q)
 \ge\int R(\beta,q)\pi_\varepsilon(\dd\beta)
 \ge r_\varepsilon.
\]
Letting $\varepsilon\downarrow0$ yields
\[
 \sup_{\beta\in\mathcal B}R(\beta,q)\ge C_J.
\]
Since this holds for every $q\in\mathcal Q$,
it proves the minimax lower bound in
Proposition~\ref{prop:minimax}.

\section{Numerical evaluation in the exponential example}
\label{app:numerical}

The calculation in Section~\ref{sec:numerical} uses the potential
\eqref{eq:example-Fourier} with $\theta_0=0$. Its Fourier multiplier
has an integrable singularity at the origin. To evaluate it, write
\[
 s=t_1^2+t_2^2+t_3^2,\qquad
 p=t_1^2t_2^2+t_1^2t_3^2+t_2^2t_3^2+t_1^2t_2^2t_3^2.
\]
Then $\prod_i(1+t_i^2)=1+s+p$, and for $t\ne0$ the multiplier splits as
\begin{equation}\label{eq:numerical-split}
 \begin{aligned}
 \frac{\pi^{3/2}e^{-s/4}}{1-\phi_{ab}(t)}
 &=\frac{\pi^{3/2}e^{-s/4}}{s}+\mathcal R(t),\\
 \mathcal R(t)
 &=\pi^{3/2}e^{-s/4}
       \left[1-\frac{p}{s(s+p)}\right].
 \end{aligned}
\end{equation}
The first term has the explicit inverse transform
\begin{equation}\label{eq:numerical-singular-part}
 h_0(u)=\frac{\sqrt\pi}{4}
            \frac{\operatorname{erf}(\norm{u})}{\norm{u}},
 \qquad h_0(0)=\frac12.
\end{equation}
Here
$\operatorname{erf}(r)=2\pi^{-1/2}\int_0^r e^{-z^2}\dd z$.
For $K=I_3$, this is exactly the leading inverse transform
$H_{ab,K}$ in \eqref{eq:tail-heat-representation}.
The formula follows by radial integration of the Gaussian multiplier.
Thus $w_{ab}=h_0+v$, where $v$ is the inverse Fourier transform
of $\mathcal R$.
The residual multiplier $\mathcal R$ is bounded near zero,
because $p=O(s^2)$, and has Gaussian decay at infinity.
Its directional limits at zero need not agree, but changing
its value at the origin does not affect its inverse Fourier
transform.

For a sampled Fourier calculation, this zero-frequency convention does
affect the quadrature. The primary convention is
$\mathcal R(0)=(4/5)\pi^{3/2}$, the angular mean of its directional
limits. Replacing it by $\pi^{3/2}$ provides a sensitivity check.
The residual multiplier is sampled on the corresponding frequency
grid and inverse-transformed to a uniform spatial grid.
The expression \eqref{eq:numerical-singular-part} is then
added at the spatial grid points.

The smoothed functions and outer expectations are computed by Fourier
convolution on that grid. For $r=1,2$, put
\[
 \chi_r(t)=\prod_{i=1}^3\frac{\Gamma(r-it_i)}{\Gamma(r)}.
\]
The multiplier for $n_r(u)=\E w_{ab}(u-L_r)$ is $\chi_r(t)$.
After taking the logarithm of the resulting positive field, the outer
expectation
$F_r(\theta)=\E\log n_r(\theta+L_r)$ uses the conjugate multiplier
$\overline{\chi_r(t)}$. Hence the calculated improvement is
$F_2(\theta)-F_1(\theta)$ on the diagonal
$\theta=(\log\beta)(1,1,1)$. This evaluates the logarithm of a
deterministic approximation to $n_r$.

The Fourier transforms were evaluated in double precision on a cube
with coordinate grid $u_j=-L+jh$, $j=0,\ldots,N-1$, where
$h=2L/N$. The primary calculation used $(N,L)=(512,32)$ and
$h=0.125$. The sampled inverse transform of $\mathcal R$ includes
the spatial normalization $h^{-3}$. All computed values of $w_{ab},n_1,n_2$
were positive before their logarithms were taken. Cubic spline
interpolation of the diagonal risk values produced the 107 evaluation
points in the displayed interval; no extrapolation was used.

Two successive domain refinements, from $(N,L)=(256,16)$ to
$(384,24)$ and then to $(512,32)$, kept $h=0.125$ fixed.
The maximum changes in the computed improvement over the evaluation
points were $1.19\times10^{-4}$ and $2.99\times10^{-5}$,
respectively. At $L=24$, decreasing $h$ from $0.125$ to
$0.09375$ changed the curve by at most $1.51\times10^{-6}$.
Changing the origin convention on the $(384,24)$ grid changed it
by at most $1.91\times10^{-5}$.

The outer expectation was also evaluated at the five values in
Table~\ref{tab:risks} using eight independently scrambled
nine-dimensional Sobol sequences, each containing $2^{16}$ points,
with seeds $20261014+k$, $k=0,\ldots,7$. The uniform coordinates
generate independent $L_2$ and $D_{ab}$: each Gamma variable of shape two
is the sum of two independent unit exponentials, and each coordinate
of $D_{ab}$ is a third unit exponential.
Setting $L_1=L_2-D_{ab}$ then gives
the beta--gamma coupling.
Cubic interpolation of the precomputed
$\log n_2$ and $\log n_1$ fields yields paired values of
\[
 \log n_2(\theta+L_2)-\log n_1(\theta+L_1).
\]
All evaluation points lay within the spatial grid, without periodic
wrapping. The largest difference between the mean of the eight
scrambled estimates and the Fourier outer expectation was
$2.07\times10^{-5}$. Each difference was at most $2.12$ estimated
standard errors, with the standard error defined as the sample standard
deviation of the eight estimates divided by $\sqrt8$.
This check uses a different outer integration procedure but shares
the same deterministic approximations to $n_1$ and $n_2$.
These checks assess numerical error empirically and do not
provide a rigorous bound on the total error.
The estimated standard errors measure variability across
the scrambled Sobol estimates, not errors in the shared
approximations to $n_1$ and $n_2$.
The numerical illustration is separate from the analytic
proof of domination.

\end{appendix}

\end{document}